\documentclass[10pt, a4paper, reqno]{amsart}
\usepackage{amscd}
\usepackage{dsfont}

\usepackage{algpseudocode,algorithm,algorithmicx}

\newcommand*\Let[2]{\State #1 $\gets$ #2}
\algrenewcommand\algorithmicrequire{\textbf{Precondition:}}
\algrenewcommand\algorithmicensure{\textbf{Postcondition:}}

\usepackage{amsmath, upgreek}
\usepackage{amstext}
\usepackage{amsthm}
\usepackage{amssymb}
\usepackage{amsopn}
\usepackage{amssymb}
\usepackage{amsxtra}
\usepackage{amsfonts}
\usepackage{fancyhdr}
\usepackage{paralist, epsfig}
\usepackage[mathcal]{euscript}
\usepackage[english]{babel}
\usepackage[latin1]{inputenc}
\usepackage{bbold}

\usepackage{tikz}
\usepackage{tkz-euclide}

\makeatletter
\g@addto@macro\normalsize{%
  \setlength\abovedisplayskip{10pt}
  \setlength\belowdisplayskip{10pt}
  \setlength\abovedisplayshortskip{5pt}
  \setlength\belowdisplayshortskip{8pt}
}

\allowdisplaybreaks

\newtheoremstyle{normal}
{5pt}
{5pt}
{\normalfont}
{}
{\bfseries}
{}
{0.4em}
{\bfseries{\thmname{#1}\thmnumber{ #2}.\thmnote{ \hspace{0.5em}(#3)\newline}}}

\newtheoremstyle{kursiv}
{5pt}
{5pt}
{\itshape}
{}
{\bfseries}
{}
{0.4em}
{\bfseries{\thmname{#1}\thmnumber{ #2}.\thmnote{ \hspace{0.5em}(#3)\newline}}}

\theoremstyle{kursiv}

\theoremstyle{normal}

\newtheorem{thm}{Theorem}
\newtheorem{ex}[thm]{Example}

\renewcommand{\epsilon}{\varepsilon}

\renewcommand{\theta}{\vartheta}

\usepackage{setspace}
\usepackage{color}

\definecolor{grey}{gray}{.3}

\begin{document}
\allowdisplaybreaks
$ $
\vspace{-40pt}

\title{On a counterexample in connection with\\the Picard-Lindel\"of theorem}

\author{Georgios Passias\hspace{0.5pt}\MakeLowercase{$^{\text{1}}$} and Sven-Ake Wegner\hspace{0.5pt}\MakeLowercase{$^{\text{2}}$}}

\renewcommand{\thefootnote}{}
\hspace{-1000pt}\footnote{\hspace{5.5pt}2010 \emph{Mathematics Subject Classification}: Primary 34A12; Secondary 65L05, 34-01.\vspace{1.6pt}}

\hspace{-1000pt}\footnote{\hspace{5.5pt}\emph{Key words and phrases}: Initial value problem, Picard-Lindel\"of theorem, Lipschitz condition, Euler method. \vspace{1.6pt}}

\hspace{-1000pt}\footnote{\hspace{0pt}$^{1}$\,University of Wuppertal, School of Mathematics and Natural Sciences, Gau\ss{}stra\ss{}e 20, 42119 Wuppertal, Germany.\vspace{1.6pt}}

\hspace{-1000pt}\footnote{\hspace{0pt}$^{2}$\,Corresponding author: Teesside University, School of Computing, Engineering \&{} Digital Technologies, Middlesbrough,\linebreak\phantom{x}\hspace{1.2pt}TS1\;3BX, United Kingdom, phone: +44\,(0)\,1642\:434\:82\:00, e-mail: s.wegner@tees.ac.uk.}

\begin{abstract}
We give an example, which demonstrates that in the situation of the Picard-Lindel\"of theorem, the Lipschitz condition on the right hand side $f(x,y)$ with respect to $y$, cannot be replaced by Lipschitz continuity in $y$ for every $x$. We show that, in our example, the classical Euler method detects only one of infinitely many solutions and we outline how the latter can be adjusted to find also other solutions numerically.
\end{abstract}

\maketitle

\vspace{-10pt}

The Picard-Lindel\"of theorem is a classical topic of undergraduate analysis. One of its most common formulations is that the initial value problem
\begin{equation}\label{CP}
y'=f(x,y),\;y(x_0)=y_0,
\end{equation}
where $x_0,\,y_0\in\mathbb{R}$, $a,\,b>0$ and $f\colon[x_0-a,x_0+a]\times[y_0-b,y_0+b]\rightarrow\mathbb{R}$ is continuous, has a unique local solution, provided that $f=f(x,y)$ satisfies a Lipschitz condition with respect to $y$, i.e., if
\begin{equation}\label{LIP}
\begin{array}{c}\vspace{4pt}
\exists\:L\geqslant0\;\forall\:x\in[x_0-a,x_0+a]\;\forall\:y_1,\:y_2\in[y_0-b,y_0+b]\colon\\
|f(x,y_1)-f(x,y_2)|\leqslant L|y_1-y_2|
\end{array}
\end{equation}
holds. From this \textquotedblleft{}rectangle version\textquotedblright{} of the theorem, variations like unique local existence for any initial condition under a local Lipschitz condition, or existence of a global solution by requiring \eqref{LIP} on a strip $[x_0-a,x_0+a]\times\mathbb{R}$, can be derived.

\medskip

In order to show that dropping the Lipschitz condition leads to loosing the uniqueness, foundational textbooks often seem to use examples like $y'=y^{2/3}$, cf.~Arnold \cite[p.~36]{A}, $y'=|y|^{1/2}$, cf.~Walter \cite[p.~15]{W}, or variants of the latter. Following this line in the classroom, a natural question, namely wether or not reordering the quantifiers in \eqref{LIP} as follows
\begin{equation}\label{LIP-2}
\begin{array}{c}\vspace{4pt}
\forall\:x\in[x_0-a,x_0+a]\;\exists\:L\geqslant0\;\forall\:y_1,\:y_2\in[y_0-b,y_0+b]\colon\\
|f(x,y_1)-f(x,y_2)|\leqslant L|y_1-y_2|
\end{array}
\end{equation}
gives a sufficient condition for uniqueness, remains untouched. Since understanding that the order of quantifiers is of utmost importance in general, and since students might in particular confuse the \emph{Lipschitz condition with respect to $y$} with \emph{Lipschitz continuity for every $x$}, we believe that a discussion of the latter question in the classroom is beneficial for students.

\medskip

The aim of this short note is firstly to give an elementary counterexample which shows that \eqref{LIP-2} is not sufficient for uniqueness. Our second aim is to show that in the situation of this counterexample the classical Euler method finds only one solution whereas a small modification of the latter then finds also other solutions. When discussing numerical techniques with students, the latter can be used to illustrate that a successful algorithm a priori needs not to reveil the full truth.

\medskip

We emphasize that there is a vast literature on the question how uniqueness could be regained by requiring \eqref{LIP-2} and additionally imposing conditions on $L=L(x)$. Results like the criteria of Nagumo, Osgood or Perron have been established more than 50 years ago and examples showing that they are sharp are well-known. We refer to Hartman \cite{H}, Agarwal, Lakshmikantham \cite{AL}, and the references therein.

\medskip

The idea behind our counterexample is to divide the first quadrant by two parabolas with vertex at zero into three \textquotedblleft{}distorted sectors\textquotedblright{}, cf.\ Fig.\ 2. The function $f$ is then defined on the top and on the bottom sector in a way that the corresponding initial value problems are integrable and their solutions stay in the respective sectors. The gap over the middle sector is filled by \textquotedblleft{}connecting opposite values with a straight line in $y$-direction\textquotedblright{}. The existence of two solutions disqualifies $f$ from satisfying a Lipschitz condition. But cutting through the graph of $f$ for fixed $x$ in $y$-direction leads to a function $f(x,\cdot)$, see Fig.\ 1, that is first constant, then linearly increasing, then again constant\,---\,and thus of course Lipschitz continuous. 

\smallskip

\begin{ex} Consider the continuous function $f\colon[0,\infty)\times[0,\infty)\rightarrow\mathbb{R}$, given by
\begin{equation}\label{FUNK}
f(x,y)=\begin{cases}
\;\hspace{34pt}x/2 &\text{if }\; 0\leqslant y\leqslant x^2/2,\\
\;x/2+5(y-x^2/2)/x &\text{if }\; x^2/2<y<x^2,\\
\;\hspace{37pt}3x &\text{if }\; x^2\leqslant{}y.
\end{cases}
\end{equation}
Then the initial value problem $y'=f(x,y)$, $y(0)=0$, has the two solutions $\varphi_1$, $\varphi_2\colon[0,\infty)\rightarrow\mathbb{R}$, $\varphi_1(x)=x^2/4$, $\varphi_2(x)=3x^2/2$, which satisfy $\varphi_1|_{[0,\varepsilon]}\not=\varphi_2|_{[0,\varepsilon]}$ for any $\varepsilon>0$. Moreover, for any $x\in[0,\infty)$ there exists $L\geqslant0$ such that for all $y_1,y_2\in[0,\infty)$ the estimate $|f(x,y_1)-f(x,y_2)|\leqslant L|y_1-y_2|$ is valid.
\end{ex}

\begin{proof}\label{PROOF} As $\varphi_1$ lies completely in the area $0\leqslant{}y\leqslant x^2/2$ and $\varphi_2$ lies completely in the area $x^2\leqslant{}y$, it is straightforward to check that $\varphi_1$ and $\varphi_2$ both solve the initial value problem. It remains to establish Lipschitz continuity for every $x$. If $x=0$, then $f(x,\cdot)\equiv0$ holds and the estimate is trivial. For $x>0$ we observe that $f(x,\cdot)$ is constant on $[0,x^2/2]$ and constant on $[x^2,\infty)$.

\begin{center}
\begin{tikzpicture}[xscale=3,yscale=0.7]
\tkzInit[xmin=0,xmax=1.2, ymin=0,ymax=20]
\draw[black, line width=0.2mm, domain=0:1/2] plot(\x,{1/2+0.005});
\draw [<->,thick, black] (0,3.5) node (yaxis) [above] {} |- (1.7,0) node (xaxis) [right] {};
\draw[black,dashed, line width=0.2mm, domain=0:3] plot({1},{\x});
\draw[black,dashed, line width=0.2mm, domain=0:1/2] plot({1/2},{\x});
\draw[black,line width=0.2mm, domain=1/2:1] plot(\x,{1/2+5*(\x-1/2)/1});
\draw[black, line width=0.2mm, domain=1:1.7] plot(\x,{3});
\draw[black, dashed, line width=0.2mm, domain=0:1/2] plot(\x,{1/2-0.02});
\draw[black,dashed, line width=0.2mm, domain=0:1] plot({\x},{3});
\draw [color=black] (-0.09,2.99) node {\footnotesize$3x$};
\draw [color=black] (-0.125,0.45) node {\footnotesize$x/2$};
\draw [color=black] (0.5,-0.36) node {\footnotesize$x^2/2$};
\draw [color=black] (1.01,-0.35) node {\footnotesize$x^2$};
\draw [color=black] (1.86,3.005) node {\footnotesize$f(x,\cdot)$};
\end{tikzpicture}
\end{center}

\vspace{-18pt}

\begin{center}
\small \textbf{Figure 1:}\:For $x>0$ the function $f(x,\cdot)$ is constant on $[0,x^2/2]$ and $[x^2,\infty)$.
\end{center}

Consequently, the estimate is again trivial, if $y_1$ and $y_2$ are both contained in one of these two intervals. If not, then we may assume w.l.o.g.\ that $x^2/2\leqslant y_1<y_2\leqslant x^2$. If $y_1<x^2/2$ we can \textquotedblleft{}move\textquotedblright{} $y_1$ rightwards to $x^2/2$ without that the values of $f(x,\cdot)$ change. Similarly, we can move $y_2>x^2$ leftwards to $x^2$. For $y_1$, $y_2\in[x^2/2,x^2]$ we finally compute $|f(x,y_1)-f(x,y_2)|=|[x/2+5(y_1-x^2/2)/x]-[x/2+5(y_2-x^2/2)/x]|=(5/x)|y_1-y_2|$ which establishes the estimate with $L=5/x$.
\end{proof}

The above can be modified as follows in order to construct an initial value problem of the form considered in \eqref{CP} with $x_0=y_0=0$, $a=b=1$. We define $f(x,y)$ as in \eqref{FUNK} if $x,y\geqslant0$, put $f(x,y)=0$ if $x<0$, and $f(x,y)=x/2$ if $x\geqslant0$ and $y<0$. Then, \eqref{CP} has the solutions $\varphi_1\colon [-1,1]\rightarrow\mathbb{R}$ and $\varphi_2\colon[-1,(2/3)^{1/2}]\rightarrow\mathbb{R}$ given by $\varphi_1(x)=x^2/4$ and $\varphi_2(x)=3x^2/2$ for $x\geqslant0$ and $\varphi_1(x)=\varphi_2(x)=0$ for $x<0$. It is straightforward to show that $f(x,\cdot)$ is Lipschitz continuous for every $x\in[-1,1]$.

\medskip

Let us now approach the Cauchy problem \eqref{CP} with $f$ as in \eqref{FUNK} from a computational point of view. Here, the Euler method yields only the solution $\varphi_1$, since $f(0,0)=0$ implies that the first part of the polygonal curve is always zero, and later the approximant stays in the area $0\leqslant y\leqslant x^2/2$ since \textquotedblleft{}the slope field keeps it away from the upper boundary of this area\textquotedblright{}. On the other hand, one can also easily reach the solution $\varphi_2$ numerically. For this, one has to modify the Euler method by pushing the approximant in the initial step into the area $x^2\leqslant y$. This can be done, e.g., via

\begin{center}
\begin{minipage}[h]{400pt}
  \begin{algorithmic}[1]
    \Statex
    \Function{PushEuler}{$f, x_0, y_0, K, h, n$}
      \State $x_1\gets x_0+h,\;y_1\gets y_0+K\cdot{}h$
      \For{$j \gets 2 \textrm{ to } n$}
      \Let{$m$}{$f(x_{j-1},y_{j-1})$}
      \State $x_{j}\gets x_{j-1}+h,\;y_{j}\gets y_{j-1}+m\cdot{}h$
      \EndFor
      \State \Return{$(x_j,\;y_j)_{j=0,\dots,n}$}
    \EndFunction
  \end{algorithmic}
 \end{minipage}
\end{center}

\smallskip

\noindent{}which leads with $K=1$ for the initial value problem of our Example to Fig.\ 2.

\medskip


\begin{center}
\begin{tikzpicture}[scale=3.8]
\tkzInit[xmin=0,xmax=1, ymin=0,ymax=1]
\fill[black!5!white]plot[domain=0:1](\x,1)--plot[domain=1:0](\x,{\x^2});
\fill[black!10!white]plot[domain=0:1](\x,{\x^2})--plot[domain=1:0](\x,{\x^2/2});
\fill[black!15!white]plot[domain=0:1](\x,0)--plot[domain=1:0](\x,{\x^2/2});
\draw[black,dashed, line width=0.2mm, domain=0:1] plot(\x,{\x^2/4});
\draw[black,line width=0.2mm, domain=0:(2/3)^(1/2)] plot(\x,{3*\x^2/2});
\draw [color=red] (0.25,0.85) node {$h=0.1$};
\draw [color=red, line width=0.2mm] plot coordinates {(0.0,0.0)(0.1,0.1)(0.2,0.16)(0.3,0.25)(0.4,0.37)(0.5,0.52)(0.6,0.7)(0.7,0.91)(0.74,1)};
\fill [color=red] (0.0,0.0)  circle[radius=0.2pt];
\fill [color=red](0.1,0.1)  circle[radius=0.2pt];
\fill [color=red](0.2,0.16) circle[radius=0.2pt];
\fill [color=red](0.3,0.25)  circle[radius=0.2pt];
\fill [color=red](0.4,0.37)  circle[radius=0.2pt];
\fill [color=red](0.5,0.52)  circle[radius=0.2pt];
\fill [color=red](0.6,0.7)  circle[radius=0.2pt];
\fill [color=red](0.7,0.91)  circle[radius=0.2pt];
\end{tikzpicture}
\begin{tikzpicture}[scale=3.8]
\tkzInit[xmin=0,xmax=1, ymin=0,ymax=1]
\fill[black!5!white]plot[domain=0:1](\x,1)--plot[domain=1:0](\x,{\x^2});
\fill[black!10!white]plot[domain=0:1](\x,{\x^2})--plot[domain=1:0](\x,{\x^2/2});
\fill[black!15!white]plot[domain=0:1](\x,0)--plot[domain=1:0](\x,{\x^2/2});
\draw[black,dashed, line width=0.2mm, domain=0:1] plot(\x,{\x^2/4});
\draw[black,line width=0.2mm, domain=0:(2/3)^(1/2)] plot(\x,{3*\x^2/2});
\draw [color=red] (0.25,0.85) node {$h=0.05$};
\draw [color=red, line width=0.2mm] plot coordinates {(0,0)(0.05,0.0500)(0.1,0.0650)(0.15,0.0875)(0.2, 0.1175)(0.25, 0.1550)(0.3,0.2000)(0.35, 0.2525)(0.4, 0.3125)(0.45, 0.3800)(0.5, 0.4550)(0.55, 0.5375)(0.6, 0.6275)(0.65, 0.7250)(0.7, 0.8300)(0.75,0.9425)(0.774,1)};
\fill [color=red] (0.0,0.0)  circle[radius=0.2pt];
\fill [color=red] (0.05,0.0500) circle[radius=0.2pt];
\fill [color=red] (0.1,0.0650) circle[radius=0.2pt];
\fill [color=red] (0.15,0.0875)  circle[radius=0.2pt];
\fill [color=red] (0.2, 0.1175) circle[radius=0.2pt];
\fill [color=red] (0.25, 0.1550) circle[radius=0.2pt];
\fill [color=red] (0.3,0.2000) circle[radius=0.2pt];
\fill [color=red] (0.35, 0.2525) circle[radius=0.2pt];
\fill [color=red] (0.4, 0.3125) circle[radius=0.2pt];
\fill [color=red] (0.45, 0.3800)  circle[radius=0.2pt];
\fill [color=red] (0.5, 0.4550) circle[radius=0.2pt];
\fill [color=red] (0.55, 0.5375) circle[radius=0.2pt];
\fill [color=red] (0.6, 0.6275) circle[radius=0.2pt];
\fill [color=red]   (0.65, 0.7250) circle[radius=0.2pt];
\fill [color=red] (0.7, 0.8300)  circle[radius=0.2pt];
\fill [color=red] (0.75,0.9425)  circle[radius=0.2pt];
\end{tikzpicture}
\begin{tikzpicture}[scale=3.8]
\tkzInit[xmin=0,xmax=1, ymin=0,ymax=1]
\fill[black!5!white]plot[domain=0:1](\x,1)--plot[domain=1:0](\x,{\x^2});
\fill[black!10!white]plot[domain=0:1](\x,{\x^2})--plot[domain=1:0](\x,{\x^2/2});
\fill[black!15!white]plot[domain=0:1](\x,0)--plot[domain=1:0](\x,{\x^2/2});
\draw[black,dashed, line width=0.2mm, domain=0:1] plot(\x,{\x^2/4});
\draw[black,line width=0.2mm, domain=0:(2/3)^(1/2)] plot(\x,{3*\x^2/2});
\draw [color=red] (0.25,0.85) node {$h=0.01$};
\draw [color=red, line width=0.2mm] plot coordinates {(0,0)(0.01,0.0100)(0.02,0.0106)(0.03,0.0115)(0.04,0.0127)(0.05,0.0142)(0.06,0.0160)(0.07,0.0181)(0.08,0.0205)(0.09,0.0232)(0.1,0.0262)(0.11,0.0295)(0.12,0.0331)(0.13,0.0370)(0.14,0.0412)(0.15,0.0457)(0.16,0.0505)(0.17,0.0556)(0.18,0.0610)(0.19,0.0667)(0.2,0.0727)(0.21,0.0790)(0.22,0.0856)(0.23,0.0925)(0.24,0.0997)(0.25,0.1072)(0.26,0.1150)(0.27,0.1231)(0.28,0.1315)(0.29,0.1402)(0.3,0.1492)(0.31,0.1585)(0.32,0.1681)(0.33,0.1780)(0.34,0.1882)(0.35,0.1987)(0.36,0.2095)(0.37,0.2206)(0.38,0.2320)(0.39,0.2437)(0.4,0.2557)(0.41,0.2680)(0.42,0.2806)(0.43,0.2935)(0.44,0.3067)(0.45,0.3202)(0.46,0.3340)(0.47,0.3481)(0.48,0.3625)(0.49,0.3772)(0.50,0.3922)(0.51,0.4075)(0.52,0.4231)(0.53,0.4390)(0.54,0.4552)(0.55,0.4717)(0.56,0.4885)(0.57,0.5056)(0.58,0.5230)(0.59,0.5407)(0.60,0.5587)(0.61,0.5770)(0.62,0.5956)(0.63,0.6145)(0.64,0.6337)(0.65,0.6532)(0.66,0.6730)(0.67,0.6931)(0.68,0.7135)(0.69,0.7342)(0.70,0.7552)(0.71,0.7765)(0.72,0.7981)(0.73,0.8200)(0.74,0.8422)(0.75,0.8647)(0.76,0.8875)(0.77,0.9106)(0.78,0.9340)(0.79,0.9577)(0.8,0.9817)(0.8075,1)};
\end{tikzpicture}

\medskip

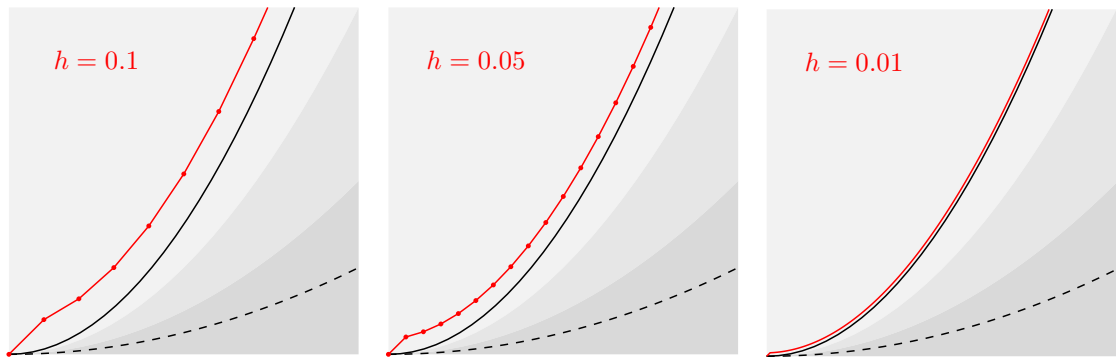
\small \textbf{Figure 2:}\:Polygonal curves (red) generated by {\sc PushEuler}, the solutions $\varphi_1$ (dashed black)\\and $\varphi_2$ (solid black). The shaded areas indicate the different formulas that define $f$.
\end{center}

\medskip

An adaption of that proof of Peano's theorem relying on the Euler method, see e.g.~\cite[p.~78]{W}, shows the following: Given a continuous right hand side $f$, then a sequence $(p_{k})_{k\in\mathbb{N}}$ of polygonal curves, corresponding to partitions of mesh size $h_{k}\rightarrow0$, and arizing from the algorithm above, contains a subsequence that converges on $[x_0,x_0+a]$ uniformly to a solution of \eqref{CP}. For $f$ as in \eqref{FUNK} and $K=1$ it is easy to show that this subsequence converges to $\varphi_2$, see \cite{P} for details.

\smallskip

We leave it to the reader to verify that the initial value problem of the example has in addition to $\varphi_1$ and $\varphi_2$ infinitely many other solutions and that (at least some of) these solutions can be approximated by our algorithm (solve the linear inhomogenous equation $y'=(5/x)y+2x$ and try, e.g., {\sc PushEuler}($f,0,0,0.07,0.1,10$)).

\medskip

{

\small 

{\sc Acknowledgements.} The authors would like to thank Lech Maligranda (Lule\aa) who drew their attention to the criteria of Nagumo, Osgood and Perron. They would also like to thank the referees and editors of College Mathematics Journal for their constructive comments.

}

\medskip

\normalsize

\bibliographystyle{amsplain}

\end{document}